\documentclass[12pt]{article}
\usepackage{amsmath}
\usepackage{amssymb}
\usepackage{latexsym}
\usepackage{epsfig}
\usepackage{graphicx}
\usepackage{mathptmx}
\usepackage[mtpcal]{}
\usepackage{mathtools}
\usepackage{physics}
\usepackage{cancel}

\begin{document}
\setlength{\baselineskip}{7mm}

\setlength{\baselineskip}{7mm}

\begin{center}
{\huge Calculus of functional centrality}         
\end{center}
\vspace*{5mm}
\begin{center}
\begin{tabular}{c}
Djemel Ziou \\
D\'epartement d'informatique  \\
 Universit\'e de Sherbrooke  \\
 Sherbrooke, Qc., Canada J1K 2R1\\
 Djemel.Ziou@usherbrooke.ca
\end{tabular}
\end{center}


\vspace*{20mm}

\begin{abstract}
In this document, we present another perspective for the calculus of optimal geometrical primitives and functions according to the centrality requirements. The shortest paths expressed in spatial and temporal domains are studied. We show the effectiveness of this formulation by providing solutions that cannot be easily accessed   by classical  formulation when using the calculus of variations.
\end{abstract}





\section{Introduction}

The calculus of variations is widely used in mathematics, physics, image processing,  and other areas of knowledge~\cite{Courant53,Benhama20,Kerouh18}.  While the optimization aims to find a stationary point of  a function defined over some domain, the calculus of variation seeks to find  stationary geometrical primitives other than the point (e.g. line, surface, volume) or a stationary  function  according to a given functional. Let us concentrate on the simplest calculus of variations in which the unknown is univariate continuously derivable curves  $u(x):[a,b] \to \mathbb{R}$ and the first-order functional to be minimized depends upon $u(x)$, its first derivative $u'(x)$, and $x$. This first-order functional is specified by  a geometrical or physics feature $F(x,u(x),u'(x))$,    such as the curve length, its curvature  or the time travelled by an object. It is given by:
\begin{equation}
J(u)= \int_{a}^{b} F(x,u(x),u'(x)) dx
\label{V1}
\end{equation}
The  feature $F(x,u(x),u'(x))$ is assumed to be  twice continuously differentiable with respect to  $x$ and the first and second variations with respect to  $u$ and  $u'$ are continuous with respect to $x$. In this case, the solution to this problem can be obtained by solving the following   Euler-Lagrange ordinary differential equation (ODE):
\begin{equation}
\frac{\partial F(x,u(x),u'(x)) }{\partial u} - \frac{ d}{ d x} \frac{\partial F(x,u(x),u'(x)) }{\partial u'} =0
\end{equation}
In order to formulate the condition so that the stationary curve to be a minimum of the functional $J(u)$, we need to establish the positive definiteness of the second variation. Straightforward computations lead to conclude that the critical  curve u(x) is a minimum if  the following inequality holds:
 \begin{equation}
\frac{\partial^2 F(x,u(x),u'(x)) }{\partial u^2} h^2 + 2 \frac{\partial^2 F(x,u(x),u'(x)) }{\partial u   \partial u'} h h'+ \frac{\partial^2 F(x,u(x),u'(x)) }{\partial u'^2} h'^2 > 0
\end{equation}
where $h(x)$ is a non-constant variation such that $h(a)=h(b)=0$. The reader  will find more about the calculus of variations in~\cite{Courant53} and in other literature. Let us now consider things from another perspective.  To this end, in this document we will limit ourselves to the case  where $F(x,u(x),u'(x))>0$. The functional in Eq.~\ref{V1} is a prior knowledge expressing some arbitrariness because it is the result of intuition not subject to experimentation. For example, it can be the length of a plane curve, which can be measured in an infinite number of ways because there are an infinite number of ways to define distance.  From this point of view, the variational problem in Eq.~\ref{V1} can be expressed in  an infinite number of ways. Let us consider the  transformation of the feature $F(x, u(x), u'(x))$ yielding $F^\alpha(x, u(x), u'(x))$, where $\alpha \in \mathbb{R}$.  For a given $x$,  $F^\alpha (x,u(x),u'(x))$ is greater or lesser than $F(x,u(x),u'(x))$ depending upon the values of $F(x, u(x), u'(x))$ and  $\alpha$. Similarly to Eq.~\ref{V1}, the first-order functional   can be written in the transformed space as: 
\begin{equation}
C_\alpha^\alpha(u)= \frac{1}{\delta} \int_{a}^{b} F^\alpha(x,u(x),u'(x)) dx
\end{equation}
where $\delta=b-a$.
The rational explanation of the exponent $\alpha$ of $C_\alpha()$ can be provided by the mean value theorem.  By considering $F()$ as an explicit or implicit function of $x$, then there exists $\eta \in ]a,b[$ such that $F^\alpha(\eta, u(\eta), u'(\eta))=\frac{1}{\delta} \int_{a}^{b} F(x,u(x),u'(x))^\alpha dx$. Instead of using $F(\eta, u(\eta), u'(\eta))$, we rename it $C(u)$.  Taking the $\alpha$-root allows us to write:
\begin{equation}
C_\alpha(u)= (\frac{1}{\delta} \int_{a}^{b} F^\alpha(x,u(x),u'(x)) dx)^{1/\alpha}
\label{V1.2}
\end{equation}
It is worth mentioning that when $\alpha=1$, the functional in Eq.~\ref{V1} is a mean value up to the interval length.  Because  $C_{\alpha}(u)$ is a limit when ($n \to +\infty$) of the H\"older mean of $n$ samples, we call it the first-order functional centrality  of H\"older. In this report, we will show that the  first-order  functional centrality of H\"older leads us to solutions that cannot be obtained by the variational formulation in Eq.~\ref{V1}. We present some of its properties in the next section. In Section 3, we derive the ODE solution to the first-order  functional centrality of H\"older. The shortest paths in both spatial and temporal domains are studied in Sections 4 and 5.  

\section{Properties}
We announce and prove the following properties  of the first-order functional centrality  of H\"older in Eq.~\ref{V1.2} under the assumption  $F(x,u(x),u'(x))>0$ and it is bounded. 


{\bf P1 :}  $\lim_{\alpha \to 0} C_\alpha(u)=e^ \frac{ \int_{a}^{b} ln F(x,u(x),u'(x))
dx}{ b-a}$.  Indeed, let us write the limit of $ln C_\alpha(u)$:
\begin{equation}
\lim_{\alpha \to 0} ln C_\alpha(u) = \lim_{\alpha \to 0} \frac{1}{\alpha} ln \frac{1}{\delta}  \int_{a}^{b} F^\alpha(x,u(x),u'(x)) dx
\label{V3}
\end{equation}
The limit has the form 0/0, so by using the L'Hopital's rule and the dominated convergence theorem, we obtain:
\begin{eqnarray}
\lim_{\alpha \to 0} ln C_\alpha(u)  & = &\lim_{\alpha \to 0}  \frac{ \int_{a}^{b} F^\alpha(x,u(x),u'(x)) ln F(x,u(x),u'(x))
dx}{ \int_{a}^{b} F^\alpha(x,u(x),u'(x)) dx} \\ 
& = & \frac{ \int_{a}^{b} ln F(x,u(x),u'(x))
dx}{\delta}
\label{V4}
\end{eqnarray}

{\bf P2:} The functional  $C_\alpha(u)$ is a increasing wrt $\alpha$. Indeed, the derivative of $C_\alpha(u)$ wrt $\alpha$
\begin{equation}
\frac{d ln C_\alpha(u) }{d \alpha} =  \frac{1}{\alpha^2} \frac{\frac{1}{\delta}  \int_{a}^{b} F^\alpha(x,u(x),u'(x)) (ln F^\alpha(x,u(x),u'(x))   -  
 ln (\frac{1}{\delta}  \int_{a}^{b} F^\alpha(x,u(x),u'(x)) dx) dx}{\frac{1}{\delta} \int_{a}^{b} F^\alpha(x,u(x),u'(x))  dx}
 \label{V6.1}
\end{equation}
Equivalently,
\begin{equation}
\frac{d ln C_\alpha(u) }{d \alpha} =  \frac{1}{\alpha^2} \frac{\frac{1}{\delta}  \int_{a}^{b} F^\alpha(x,u(x),u'(x)) ln \frac{F^\alpha(x,u(x),u'(x))  }{\frac{1}{\delta}  \int_{a}^{b} F^\alpha(x,u(x),u'(x)) dx}dx}{\frac{1}{\delta} \int_{a}^{b} F^\alpha(x,u(x),u'(x))  dx}
\label{V6.2}
\end{equation}
Since $ln r \ge  1-1/r$ for all positive $x$, where the equality is true only if $r=1$. In our case, $r=F^\alpha(x,u(x),u'(x)) /$ $ \frac{1}{\delta}  \int_{a}^{b} F^\alpha(x,u(x),u'(x)) dx$ cannot be one because $F^\alpha(x,u(x),u'(x))$ cannot be a constant. Hence,
\begin{align}
\frac{d ln C_\alpha(u) }{d \alpha} >   \frac{1}{\alpha^2} \frac{\frac{1}{\delta}  \int_{a}^{b} F^\alpha(x,u(x),u'(x)) (1- \frac{\frac{1}{\delta}  \int_{a}^{b} F^\alpha(x,u(x),u'(x)) dx}{F^\alpha(x,u(x),u'(x))  })dx}{\frac{1}{\delta} \int_{a}^{b} F^\alpha(x,u(x),u'(x))  dx}
\end{align}
It can also be rewritten as follows:
\begin{align}
\frac{d ln C_\alpha(u) }{d \alpha} >   \frac{1}{\alpha^2} \frac{\frac{1}{\delta}  \int_{a}^{b} F^\alpha(x,u(x),u'(x)) dx - \frac{1}{\delta}  \int_{a}^{b} F^\alpha(x,u(x),u'(x)) dx \frac{1}{\delta}  \int_{a}^{b} \frac{F^\alpha(x,u(x),u'(x)) }{F^\alpha(x,u(x),u'(x))  }dx }{\frac{1}{\delta} \int_{a}^{b} F^\alpha(x,u(x),u'(x))  dx}
\end{align}
Calculating the leftmost integral allows us to write:
\begin{align}
\frac{d ln C_\alpha(u) }{d \alpha} >   \frac{1}{\alpha^2} \frac{\frac{1}{\delta}  \int_{a}^{b} F^\alpha(x,u(x),u'(x)) dx - \frac{1}{\delta}  \int_{a}^{b} F^\alpha(x,u(x),u'(x)) dx }{\frac{1}{\delta} \int_{a}^{b} F^\alpha(x,u(x),u'(x))  dx} = 0
\end{align}
The derivative of $ln C_\alpha(u)$ and $C_\alpha(u)$ are both positives  meaning $C_\alpha(u)$ is increasing in $\alpha$.

 
 {\bf P3:}  $C_{-\infty}(u) \le C_\alpha(u)  \le  C_{+\infty}(u)$. This is a corollary of $P2$. 
 
 {\bf P4:}  $C_{-\infty}(u) = \inf F(x,u(x),u'(x))$ and  $C_{+\infty}(u) = \sup F(x,u(x),u'(x))$.  Let  the sequence $\alpha_1,\cdots, \alpha_n$ such that   $\alpha_1 < \cdots < \alpha_n$, so according to  P2 and Eq.~\ref{V1.2}, $C_{\alpha_1}(u) < \cdots <  C_{\alpha_n}(u)$. Thanks to the mean value theorem,    the last inequality can be rewritten $C_{\alpha_1}(u)=F(\eta_1, u(\eta_1), u'(\eta_1)) < \cdots <  C_{\alpha_n}(u)=F(\eta_n, u(\eta_n), u'(\eta_n)) $. The sequence  is built such that  $\lim_{n \to \infty} \alpha_n = \alpha$. If   $\lim_{\alpha \to +\infty} C_{\alpha} (u) = F(\eta, u(\eta), u'(\eta)) $ is finite, then it is the upper bound because $C_\alpha(u)$ is increasing in $\alpha$. The same reasoning is used for the $\lim_{\alpha \to -\infty} C_\alpha (u)$. 
  
The properties P2, P3 and P4 ensure that the $C_\alpha(u)()$ is non decreasing in $\alpha$ and it is bounded, i.e.   $C_\alpha(u) \in [\inf F(), \sup F()]$ for a given $\alpha$.     The parameter $\alpha$ implements a transform of the features $F( )$ to a new space in which the first-order functional centrality of H\"older  is defined, i.e. Eq.~\ref{V1.2}. This transform leads to specifying the contribution of the point $x$ to H\"older's first-order functional centrality using the parameter $\alpha$. In other words, the relevance of features is specified by the parameter $\alpha$, meaning that it is a trick for feature (data) selection~\cite{Ziou23}. Taken together these properties indicate that  the functional $C_\alpha(u)$ is a mean operator. The discrete case (i.e., counting measure) of  this functional is the H\"older mean. By the analogy of the discrete case, we call it  the functional arithmetic mean when $\alpha=1$,  the functional harmonic mean when $\alpha=-1$,  and  according to P1 the functional geometric mean when $\alpha=0$.
  
\section{Solutions of the First order functional centrality of H\"older}
In a space of curves, we want to find  a smooth,  continuously derivable curve $u(x)$ defined over $[a,b]$  minimizing $C_\alpha(u)$ in Eq.~\ref{V1.2}. The curve is assumed to be known in two points, $P=(a,u(a))$ and $Q=(b,u(b))$, where the values $a$ and $b$ are  the data of the problem to be studied. In order to find the stationary  curve $u(x)$,  we vary $u(x)$  by a small amount of $h(x)$ (i.e. $u_\epsilon(x)= u(x) + \epsilon h(x)$ and $u'_\epsilon(x)= u'(x) + \epsilon h'(x)$), and calculate  the  first variation wrt $\epsilon$. Note that, all variations are "glued" at the two points $P$ and $Q$; that is the real function $h(x)$ must verifies $h(a)=h(b)=0$.  According to P1 and Eq.~\ref{V1.2}, the  two cases $\alpha \neq 0$ and   $\alpha=0$ are  considered separately. 


\subsection{Case $\alpha \neq 0$}
The  first variation wrt $\epsilon$ is:
\begin{equation}
\frac{d C_\alpha(u_\epsilon) }{d \epsilon} = \frac{d}{d \epsilon}  (\frac{1}{\delta}\int_a^b F(x, u_\epsilon, u_\epsilon^{'})^\alpha dx)^{1/\alpha}
\end{equation}
Applying derivative rules, we obtain:
\begin{align*}
\frac{d C_\alpha(u_\epsilon)}{d \epsilon}  &=   \frac{1}{\alpha} (\frac{1}{\delta}\int_a^b F(x, u_\epsilon, u_\epsilon^{'})^\alpha dx)^{1/\alpha -1} \\ & \frac{1}{\delta}\int_a^b \alpha F(x, u_\epsilon, u_\epsilon^{'})^{\alpha-1} \frac{d F(x, u_\epsilon, u_\epsilon^{'})}{d \epsilon}  dx =0
\end{align*}
  Let us recall that $F()>0$, and therefore  $
(\int_a^b F(x, u_\epsilon, u_\epsilon^{'})^\alpha dx)^{1/\alpha -1} \neq 0$. So, the term is:
\begin{equation}
\int_a^b  F(x, u_\epsilon, u_\epsilon^{'})^{\alpha-1} \frac{d F(x, u_\epsilon, u_\epsilon^{'})}{d \epsilon}  dx =0
\label{ode1}
\end{equation}
The total derivative wrt to $\epsilon$ is:
\begin{equation}
 \frac{d}{d \epsilon} F(x, u_\epsilon, u_\epsilon^{'})  = \frac{ \partial F(x, u_\epsilon, u_\epsilon^{'})}{\partial x} \frac{ \partial x}{\partial \epsilon}   + \frac{ \partial F(x, u_\epsilon, u_\epsilon^{'})}{\partial u_\epsilon} \frac{ \partial u_\epsilon}{\partial \epsilon} + \frac{ \partial F(x, u_\epsilon, u_\epsilon^{'})}{\partial u_\epsilon^{'}} \frac{ \partial u_\epsilon^{'}}{\partial \epsilon}
\end{equation}
Because  $\frac{ \partial x}{\partial \epsilon} =0$, by substituting this equation in Eq.~\ref{ode1} and setting $\epsilon=0$, we obtain:
\begin{equation}
\int_a^b  F(x, u, u^{'})^{\alpha-1} (\frac{ \partial F(x, u, u^{'})}{\partial u} h(x) + \frac{ \partial F(x, u, u^{'})}{\partial u^{'}} h^{'}(x)) dx =0
\end{equation}
The  integration by part of the second term leads to:
\begin{equation}
\int_a^b  h(x) (F(x,u, u^{'})^{\alpha-1} \frac{ \partial F(x, u, u^{'})}{\partial u}  -  \frac{d}{d x} F(x, u, u^{'})^{\alpha-1} \frac{ \partial F(x, u, u^{'})}{\partial u^{'}})  dx =0
\end{equation}
The fundamental lemma of  the calculus of variations~\cite{Courant53} ensures that the stationary curve $u(x)$ is a solution to  the ODE:
\begin{equation}
F(x, u, u^{'})^{\alpha-1} \frac{ \partial F(x, u, u^{'})}{\partial u}  - \frac{d}{d x} F(x, u, u^{'})^{\alpha-1} \frac{ \partial F(x, u, u^{'})}{\partial u^{'}} =0
\end{equation}
This ODE can be rewritten as:
\begin{equation}
F(x, u, u^{'})^{\alpha-2} (F(x, u, u^{'}) (\frac{ \partial F(x, u, u^{'})}{\partial u} -
 \frac{d}{d x}  \frac{ \partial F(x, u, u^{'})}{\partial u^{'}}) -  
(\alpha-1) \frac{ \partial F(x, u, u^{'})}{\partial u^{'}} \frac{d}{d x} F(x, u, u^{'}) ) =0
\end{equation}
Because $F(x, u, u^{'})^{\alpha-2} \ne 0$, the final ODE is written as follows:
\begin{equation}
F(x, u, u^{'}) (\frac{ \partial F(x, u, u^{'})}{\partial u}  -  \frac{d}{d x}  \frac{ \partial F(x, u, u^{'})}{\partial u^{'}})-
(\alpha-1) \frac{ \partial F(x, u, u^{'})}{\partial u^{'}} \frac{d}{d x} F(x, u, u^{'})  =0
\label{odeFinal}
\end{equation}
Note that if $\alpha=1$,  we have the well-known Euler-Lagrange ordinary differential equation (ODE). 
 The presence of  the term $\frac{ \partial F}{\partial u^{'}} \frac{d}{d x} F(x, u, u^{'}) $ when  $\alpha \neq 1$  could leads to other solutions. To assert that the stationary curve $u(x)$  solution to Eq.~\ref{odeFinal} is a minimum, we need to test the positiveness of the  second variation wrt $\epsilon$. Indeed, the curve $u(x)$ is a minimum if $C_\alpha(u) < C_\alpha(u+\epsilon h)$. By using Maclaurin expansion wrt to $\epsilon$, we write:
\begin{equation}
f(\epsilon) = C_\alpha(u+\epsilon h) = C_\alpha(u)+ \epsilon f'(0) + 0.5 \epsilon^2 f''(0) + O_3(u)
\end{equation}
If $u(x)$ is a stationary curve then  the first variation $f'(0)=0$. Consequently, the functional $C_\alpha(u+\epsilon h) > C_\alpha(u)$ if  $f''(0)$ is positive. At $\epsilon=0$, the first variation can be rewritten as:
\begin{equation}
\frac{ d C_\alpha(u_\epsilon)}{d \epsilon}|_{\epsilon=0}= f'(0)=
\frac{1}{\alpha} (\frac{1}{\delta} \int_{a}^{b} F^\alpha(x, u, u^{'}) dx)^{\frac{1}{\alpha}-1}  (\frac{1}{\delta} \int_{a}^{b} \frac{d  F^\alpha(x, u_\epsilon, u_\epsilon^{'})}{d \epsilon}|_{\epsilon=0} dx) 
\end{equation}
 The second variation is:
\begin{equation}
\frac{ d^2 C_\alpha(u_\epsilon)}{d \epsilon^2}|_{\epsilon=0}  = f''(0)= \begin{array}{ll}
& \frac{1}{\alpha} \frac{1}{\delta^2}(\frac{1}{\delta} \int_{a}^{b} F^\alpha(x, u, u^{'}) dx)^{\frac{1}{\alpha}-2}  ( \int_{a}^{b} F^\alpha(x, u, u^{'}) dx  \int_{a}^{b} \frac{d^2  F^\alpha(x, u_\epsilon, u_\epsilon^{'})}{d \epsilon^2}|_{\epsilon=0} dx +  \\     &  ~~~~(\frac{1}{\alpha}-1) ( \int_{a}^{b} \frac{ dF^\alpha(x, u_\epsilon, u_\epsilon^{'})}{d \epsilon}|_{\epsilon=0} dx)^2 )  
\end{array} 
\end{equation}
The stationary curve is a minimum if $\frac{ d^2 C_\alpha(u_\epsilon)}{d \epsilon^2}|_{\epsilon=0}>0$. Knowing that  both $\frac{1}{\delta}$ and $ \int_{a}^{b} F^\alpha(x, u, u^{'}) dx $ are positive,  and $\int_{a}^{b} \frac{ dF^\alpha(x, u_\epsilon, u_\epsilon^{'})}{d \epsilon}|_{\epsilon=0} dx=0$,  it follows that the stationary curve is a minimum if:
\begin{equation}
\frac{1}{\alpha} \int_{a}^{b} \frac{d^2  F^\alpha (x, u_\epsilon, u_\epsilon^{'})}{d \epsilon^2}|_{\epsilon=0} dx  >  0
\label{TestFinalN0}
\end{equation}
The inequality is fulfilled when  $\alpha$ and $\int_{a}^{b} \frac{d^2  F^\alpha (x, u_\epsilon, u_\epsilon^{'})}{d \epsilon^2}|_{\epsilon=0} dx$  have the same sign.
Two important cases that should be studied   are the explicit  independence of $F(x,u,u^{'})$ on $u(x)$ and $x$. For the former case, substituting $\frac{\partial F(x, u, u^{'})}{\partial u}=0$ in  Eq.~\ref{odeFinal} gives $ F(x, u, u^{'})  \frac{d}{d x}  \frac{ \partial F(x, u, u^{'})}{\partial u^{'}}+
(\alpha-1) \frac{ \partial F(x, u, u^{'})}{\partial u^{'}} \frac{d}{d x} F(x, u, u^{'})  =0$. An interpretation in the general case seems to be not easy to draw. However, when $\alpha=1$,  then  $\frac{d}{d x}  \frac{ \partial F(x, u, u^{'})}{\partial u^{'}}=0$ and therefore $  \frac{ \partial F(x, u, u^{'})}{\partial u^{'}}$ is a constant. This  outcome is known as the ignorable coordinates in the calculus of variations. When $\alpha=2$,  Eq.~\ref{odeFinal} is rewritten as $\frac{d }{d x}F(x, u, u^{'})   \frac{ \partial F}{\partial u^{'}}=0$, then $F(x, u, u^{'})   \frac{ \partial F}{\partial u^{'}}$ is a constant.  For the  last $ \frac{ \partial F(x, u, u^{'})}{\partial x} =0$, straightforward manipulations leads to the ODE  $F(x, u, u^{'}) \frac{d}{d x} (u' \frac{ \partial F(x, u, u^{'})}{\partial u^{'}} - F(x, u, u^{'})) +  (\alpha-1) u' \frac{ \partial F(x, u, u^{'})}{\partial u^{'}} \frac{d}{d x} F(x, u, u^{'})  =0$. When $\alpha=1$,  we obtain the well-known Beltrami identity; i.e. $u' \frac{ \partial F(x, u, u^{'})}{\partial u^{'}} - F(x, u, u^{'})$ is a constant. When $\alpha=2$, we have $\frac{d}{d x} F(x, u, u^{'})  u' \frac{ \partial F(x, u, u^{'})}{\partial u^{'}} = F(x, u, u^{'}) \frac{d}{d x} F(x, u, u^{'})$ and hence $  u' \frac{ \partial F(x, u, u^{'})}{\partial u^{'}} - \frac{1}{2} F(x, u, u^{'})$ is a constant. The case $\alpha=2$ generalizes the Beltrami identity. The study of other values of $\alpha$ could lead to other outcomes.

\subsection{Case $\alpha=0$}
We will now give the first   and second variations in the case where $\alpha=0$. The first order variation is:
\begin{align}
\frac{d C_\alpha(u)}{d \epsilon}|_{\epsilon=0} =C_\alpha(u) \int_{a}^b (\frac{\partial ln F(x, u, u^{'})}{\partial u} h + \frac{\partial ln F(x, u, u^{'})}{\partial u'} h')dx =0
\label{odeFinala01}
\end{align}
The corresponding ODE is:
\begin{equation}
\frac{ \partial ln F(x, u, u^{'})}{\partial u}  -  \frac{d}{d x}  \frac{ \partial ln F(x, u, u^{'})}{\partial u^{'}} =0
\label{odeFinala00}
\end{equation}
If $u(x)$ is an  ignorable coordinate, then   $\frac{ \partial ln F(x, u, u^{'})}{\partial u^{'}}$ is a constant.  This is another outcome generalizing the Beltrami identity. When $F(x, u, u^{'})$ is explicitly independent upon $x$, then we deduce that $ln F(x, u, u^{'})-u' \frac{ \partial ln  F(x, u, u^{'})}{\partial u^{'}}$ is a constant.  Straightforward manipulations lead to writing the second variation:
\begin{equation}
\begin{aligned}
\frac{d^2 C_\alpha(u)}{d \epsilon^2}|_{\epsilon=0} &=C_\alpha(u)((\int_{a}^b  (\frac{\partial ln F(x, u, u^{'})}{\partial u} h + \frac{\partial ln F(x, u, u^{'})}{\partial u'} h')dx)^2+ \int_{a}^b  (\frac{1}{F(x, u, u^{'})} (\frac{\partial^2 F(x, u, u^{'})}{\partial u^2} h^2+\\  & 2 \frac{\partial^2 F(x, u, u^{'})}{\partial u \partial u'} h h'+ \frac{\partial^2 F(x, u, u^{'})}{\partial u'^2} h'^2) -  \frac{1}{F^2(x, u, u^{'})}(\frac{\partial F(x, u, u^{'})}{\partial u} h + \frac{\partial F(x, u, u^{'})}{\partial u'} h')^2)dx)
\label{odeFinala02}
\end{aligned}
\end{equation}
The second variation at the stationary curve is obtained by substituting Eq.~\ref{odeFinala01} in Eq.~\ref{odeFinala02}:
\begin{equation}
\begin{aligned}
\frac{d^2 C_\alpha(u)}{d \epsilon^2}|_{\epsilon=0} = & C_\alpha(u) \int_{a}^b  (\frac{1}{F(x, u, u^{'})} (\frac{\partial^2 F(x, u, u^{'})}{\partial u^2} h^2+ 2 \frac{\partial^2 F(x, u, u^{'})}{\partial u \partial u'} h h'+ \frac{\partial^2 F(x, u, u^{'})}{\partial u'^2} h'^2) \\  & -  \frac{1}{F^2(x, u, u^{'})}(\frac{\partial F(x, u, u^{'})}{\partial u} h + \frac{\partial F(x, u, u^{'})}{\partial u'} h')^2)dx
\label{ode0st}
\end{aligned}
\end{equation}
Knowing that $C_\alpha(u)>0$, the second variation test leads to concluding that the stationary curve is a minimum if:
\begin{equation}
\begin{aligned}
 \int_{a}^b  & \frac{1}{F(x, u, u^{'})} (\frac{\partial^2 F(x, u, u^{'})}{\partial u^2} h^2+ 2 \frac{\partial^2 F(x, u, u^{'})}{\partial u \partial u'} h h'+ \frac{\partial^2 F(x, u, u^{'})}{\partial u'^2} h'^2) dx > \\ &  \int_{a}^b \frac{1}{F^2(x, u, u^{'})}(\frac{\partial F(x, u, u^{'})}{\partial u} h + \frac{\partial F(x, u, u^{'})}{\partial u'} h')^2)dx
\end{aligned}
\end{equation}
The RHS is non-negative, so we can write  the second variation test as follows:
\begin{equation}
\begin{aligned}
& \int_a^b \frac{1}{F(x, u, u^{'})}(\frac{\partial^2 F(x, u, u^{'})}{\partial u^2} h^2 + 2 \frac{\partial^2 F(x, u, u^{'})}{\partial u \partial u'} h h'+ \frac{\partial^2 F(x, u, u^{'})}{\partial u'^2} h'^2) > 0 ~~\mbox{if}~~\alpha=0
\label{TestFinal0}
\end{aligned}
\end{equation}

\section{Shortest path problem}
\begin{figure}
\centering
\includegraphics[scale=0.4]{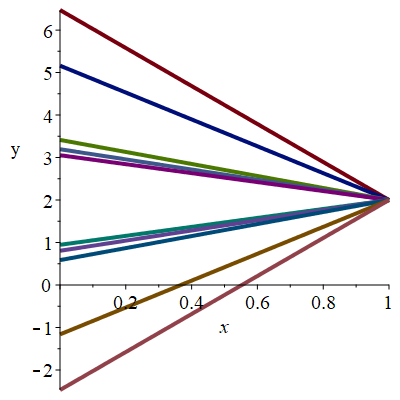} 
\caption{Real solution to the functional centrality problem.  The lines crossing the y-axis at position greater than  (resp. less than) 2 have a negative (resp. positive) slope. }
\label{LineReal}
\end{figure}

We would like to find the curve of the shortest path between the two points $P$ and $Q$ in the spatial domain. We can approximate the length of a plane curve defined on $[a,b]$ by subdividing it into infinitesimal linear pieces, measuring the length of each, and adding up  all the lengths.  The infinitesimal  piece length between $(x, u(x))$ and $(x + dx, u(x+dx))$ is the length of the hypotenuse of the right triangle having  sides $dx$ and $du(x)$, i.e. $\sqrt{dx^2+(u(x+dx)-u(x))^2} = dx \sqrt{1+\frac{(u(x+dx)-u(x))^2}{dx^2}}$. When $dx$ is too small,  the infinitesimal length is $\sqrt{1+u(x)'^2} dx$.  Let us define the feature $F(x,u, u')=\sqrt{1+u(x)'^2}$. It should be noted that $u(x)$ and $x$ are ignorable coordinates, i.e.  $\frac{ \partial F(x, u, u^{'})}{\partial u}  =0$ and  $\frac{ \partial F(x, u, u^{'})}{\partial x}  =0$. The other terms involved  in Eq.~\ref{odeFinal} are $\frac{ \partial F(x, u, u^{'})}{\partial u^{'}} = \frac{u'}{F}$,  $ \frac{d F(x, u, u^{'})}{d x}  = \frac{ u' u''}{F} $, and $ \frac{ d }{d x}  \frac{ \partial F(x, u, u^{'})}{\partial u^{'}} =\frac{F(x, u, u^{'})u'' - \frac{u'^2 u''}{F(x, u, u^{'})}}{F^2(x, u, u^{'})}$. The ODE in  Eq.~\ref{odeFinal} is rewritten as:
\begin{equation}
-F(x, u, u^{'}) \frac{F(x, u, u^{'})u'' - \frac{u'^2 u''}{F(x, u, u^{'})}}{F^2(x, u, u^{'})}-(\alpha-1) \frac{u'}{F(x, u, u^{'})} \frac{ u' u''}{F(x, u, u^{'})}=0
\end{equation}
\begin{equation}
-F(x, u, u^{'}) (F(x, u, u^{'})u'' - \frac{u'^2 u''}{F(x, u, u^{'})})-(\alpha-1) u'^2  u''=0
\end{equation}
\begin{equation}
u''(x) (1+(\alpha-1) u'^2(x))=0
\end{equation}
There are two cases $1+(\alpha-1) u'^2=0$ and $u''=0$. The solutions to these ODE are provided in Table~\ref{ODE_SP} as function of $\alpha$. The solution to the ODE $u''=0$ is a straight line  in real space completely defined by the two points, P and Q through which it passes, i.e $y=sx+d$. Although the straight line does not depend on $\alpha$, the second variation test does. It is a minimum except for the cases when the slope $s \notin ]-1,1[$ and  $\alpha < 1$.  When $\alpha < 0$  or  $\alpha \in ]0,1[$,  the  slopes of the two straight lines  solutions  to the ODE $1+(\alpha-1) u'^2=0$ are real and  $\alpha$ dependant. However, these solutions are minima when $\alpha < 0$ and maxima when $\alpha \in ]0,1[$. Let us make explicit the case of the functional harmonic mean $\alpha=-1$ to indicate that this property remains valid. For the functional geometric mean $\alpha=0$, the solutions to the ODE are two straight lines that do not depend on $\alpha$ and whose second variation test is inconclusive. The functional $F(x,u,u')$ is not real  when $\alpha \in ]1,2[$. For the case when $\alpha \ge 2$, there are two lines in the complex plane which are solutions to the ODE and which are minima. To illustrate, figure~\ref{LineReal} presents a bundle of lines passing through $P=(1,2)$ in the case where $\alpha \in ]0,1[$. The slopes of the lines crossing the $y$-axis at position greater than  (resp. less than) two are $- 1/\sqrt{1-\alpha}$ (resp. $1/\sqrt{1-\alpha}$)  for $\alpha= 0.1, 0.3,0.5, 0.9, 0.95$. Finally, it should be noted that  the functional arithmetic mean $\alpha=1$ leads to only one ODE, $u''=0$. It is  the case that is studied  in the calculus of variations. All the other cases are uncommon.

 \begin{table}
 \begin{tabular}{|c|c|c|c|c|c|} \hline
$\alpha$ & ODE &  $F(x,u,u')$  & 2nd variation test & $u(x)$ & space   \\ \hline \hline
    $>2$& $1+(\alpha-1) u'^2(x)=0$ & $\sqrt{1+\frac{1}{\alpha-1}}$  & $ \frac{1}{\delta} \frac{(\alpha-1)^3}{\alpha^2(\alpha-2)} (\frac{\alpha(\alpha-2)}{(\alpha-1)^2})^{\alpha/2}	   \int_a^b h'^2(x) dx >0$   & $ \frac{ \mp i x }{\sqrt{\alpha-1}} +d$ & Complex   \\ 
  &  &  &  minimum & &   \\ 
     $ $& $u''(x)=0$ & $\sqrt{1+s^2}$  & $ \frac{1}{\delta}  (s^2+1)^{\alpha/2-2}((\alpha-1)s^2+1)	  \int_a^b h'^2(x) dx >0$   & $ sx +d$ & Real   \\ 
  &  &  &  minimum  & &   \\ \hline 
    $2$& $1+ u'^2(x)=0$ & $0$  & $ \frac{1}{\delta}    \int_a^b h'^2(x) dx >0$   & $  \mp i x  +d$ & Complex   \\ 
  &  &  &  minimum & &   \\ 
     $ $& $u''(x)=0$ & $\sqrt{1+s^2}$  & $ \frac{1}{\delta  \sqrt{s^2+1}} 	  \int_a^b h'^2(x) dx >0$   & $ sx +d$ & Real   \\ 
  &  &  &  minimum  & &   \\ \hline 
 $ ]1,2[$& $1+(\alpha-1) u'^2(x)=0$ & Undefined  & Undefined   & $ \frac{\mp i x}{\sqrt{\alpha-1}}  +d$ & Complex  \\  
   $ $& $u''(x)=0$ & $\sqrt{1+s^2}$  & $ \frac{1}{\delta} (s^2+1)^{\alpha/2-2}((\alpha-1)s^2+1)	  \int_a^b h'^2(x) dx >0$   & $ sx +d$ & Real   \\ 
  &  &  &  minimum  & &   \\ \hline 
 1& $u''(x)=0$ & $\sqrt{1+s^2}$  & $\frac{1}{\delta (1+s^2)^{3/2}} \int h'^2(x) dx >0 $  & $s x +d$ & Real   \\ 
  &  &  &  minimum  & &   \\ \hline
 $ ]0,1[$& $1+(\alpha-1) u'^2(x)=0$ & $\sqrt{1+\frac{1}{1-\alpha}}$  & $ \frac{\alpha}{\delta} (\frac{(\alpha-1)^2+1}{(\alpha-1)^2})^{\alpha/2} \frac{(\alpha-1)^3}{((\alpha-1)^2+1)^2}	  \int_a^b h'^2(x) dx <0$   & $ \frac{\mp x}{\sqrt{1-\alpha}} +d$ & Real   \\ 
  &  &  &  maximum & &   \\ 
   $ $& $u''(x)=0$ & $\sqrt{1+s^2}$  & $ \frac{1}{\delta} (s^2+1)^{\alpha/2-2}((\alpha-1)s^2+1)	  \int_a^b h'^2(x) dx >0$   & $ sx +d$ & Real   \\ 
  &  &  &  minimum if $s^2 < 1$ & &   \\ \hline 
  0&$u''(x)=0$ & $\sqrt{1+s^2}$  & $\frac{1}{\delta} \frac{1-s^2}{(1+s^2)^2} \int h'^2(x) dx >0 $  & $s x +d$ & Real   \\ 
  &  &  &  minimum  if $s \in ]-1,1[$  & &   \\
   &  &  &   maximum if $s  \in ]-\infty, -1[ \cup ]1, +\infty[$  & &   \\ 
    & $u'(x)= \mp 1$ & $\sqrt{2}$  & 0  & $\mp x +d$ & Real   \\ 
  &  &  &  inconclusive & &     \\ \hline 
  $<0$& $1+(\alpha-1) u'^2(x)=0$ & $\sqrt{1+\frac{1}{1-\alpha}}$  & $ \frac{\alpha}{\delta} (\frac{(\alpha-1)^2+1}{(\alpha-1)^2})^{\alpha/2} \frac{(\alpha-1)^3}{((\alpha-1)^2+1)^2}	  \int_a^b h'^2(x) dx >0$   & $ \frac{\mp x}{\sqrt{1-\alpha}} +d$ & Real   \\ 
  &  &  &  minimum & &   \\ 
     $ $& $u''(x)=0$ & $\sqrt{1+s^2}$  & $ \frac{1}{\delta} (s^2+1)^{\alpha/2-2}((\alpha-1)s^2+1)	  \int_a^b h'^2(x) dx >0$   & $ sx +d$ & Real   \\ 
  &  &  &  minimum if $s \in ]-1,1[$ & &   \\ \hline 
\end{tabular}
\caption{Solutions to the ODE $u'' (1+(\alpha-1) u'^2)=0$ expressed as function of $\alpha$.}
\label{ODE_SP}
 \end{table}



\section{Fastest path problem}
The Brachistochrone is an old problem known  in the calculus of variations~\cite{Courant53}. It describes a curve that carries a particle under gravity  from one height to another in minimal time. Other variations
of this problem have included the effects of friction,
the motion of a disc on a hemisphere, the motion of a cyclist in a velodrome, and even the quantum Brachistochrone problem~\cite{Benhama20}. 

\begin{figure}[h]
\center
\includegraphics[width=5cm]{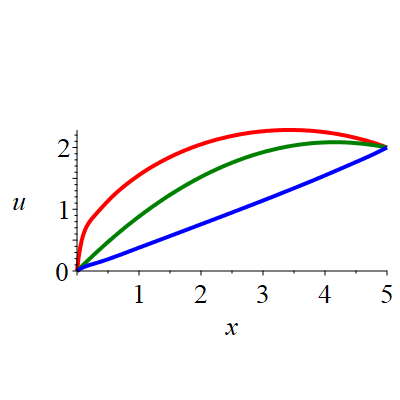}{a}
\includegraphics[width=5cm]{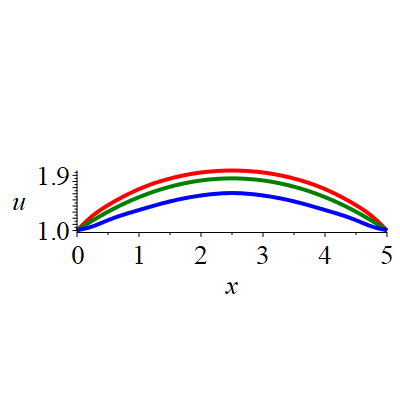}{b}
\caption{Numerical solution to Brachistochrone with the initial condition: a) u(0)=0 and u(5)=2, b) u(0)=1 and u(5)=1,  $\alpha=1$ (Red), $\alpha=2$ (Green),  c) $\alpha=8$ (Blue).} 
\label{SolTime1}
\end{figure}

Suppose that the particle is released  from $P$ at $t=0$, and then follows a curve $u(x)$ which reaches $Q$, so that $u(b)-u(a)$ is the height lost, and $b-a$ is the horizontal distance traversed. The associated feature is $F(x, u, u^{'})= \sqrt{(1+u'^2(x))/u(x)} $.  It should be noted that $\frac{ \partial F(x, u, u^{'})}{\partial x}  =0$. The other derivatives in Eq.~\ref{odeFinal} are  $ \frac{\partial F(x, u, u^{'})}{\partial u} =-\frac{\sqrt{1+u'^2}}{2 u^{3/2}}$, $\frac{\partial F(x, u, u^{'})}{\partial u'} =\frac{u'}{\sqrt{u(1+u'^2)}} $, 
 $\frac{d }{d x} \frac{\partial F(x, u, u^{'})}{\partial u'} =\frac{1}{\sqrt{u(1+u'^2)}} (\frac{u''}{1+u'^2}-\frac{u'^2}{2 u})$, and $\frac{d F(x, u, u^{'}) }{d x}  =\frac{u'u''}{\sqrt{u(1+u'^2)}}-\frac{u' \sqrt{1+u'^2}}{2u^{3/2}}$, 
$F(x, u, u^{'})( \frac{\partial F(x, u, u^{'})}{\partial u} - \frac{d }{d x} \frac{\partial F(x, u, u^{'})}{\partial u'})=-\frac{1+u'^2}{2 u^{2}}-\frac{1}{u} (\frac{u''}{1+u'^2}-\frac{u'^2}{2 u}) = -\frac{1}{2 u^{2}}- \frac{u''}{u(1+u'^2)}$, 
$(\alpha-1)\frac{\partial F(x, u, u^{'})}{\partial u'}\frac{d F(x, u, u^{'}) }{d x} =(\alpha-1) \frac{u'^2}{u} (\frac{u''}{1+u'^2}-\frac{1}{2 u})$.
The Euler-Lagrange ODE in Eq.~\ref{odeFinal} becomes:
\begin{equation}
-\frac{1}{2 u}-\frac{u''}{1+u'^2}-(\alpha-1)u'^2(\frac{u''}{1+u'^2}-\frac{1}{2 u})  =0
\end{equation}
The  solution in the calculus of variations (i.e. $\alpha=1$) is the  cycloid with the parametric equations  $x=r (\theta - sin \theta)$ and  $u(x) = r (1 - cos \theta)$.  It can be displayed by the path traced out by a point placed on a rolling wheel of radius $r$ when it rotates by an angle $\theta$. The general solution to the  ODE when $\alpha \neq 1$ is challenging. For example, it is easy to check that $u(x)= \frac{1}{\sqrt{\alpha -1}} x$ is a stationary  straight line.  According to the second variation test in Eq.~\ref{TestFinalN0}, it is a minimum  when $\alpha=2$ and  $x>0$, in this case the second variation $\frac{1}{\delta} \int_{a}^{b} ((xh'(x)-h(x))^2 +h(x))/x^3 dx>0$. We  solve the ODE numerically, using the midpoint method of Maple for three values of $\alpha$ and two different initial conditions (see  Fig.~\ref{SolTime1}). Note that the solution resembles  the Trochoid at different distances from the center of a circle of some radius. 

\begin{figure}[h]
\center
\includegraphics[width=4cm]{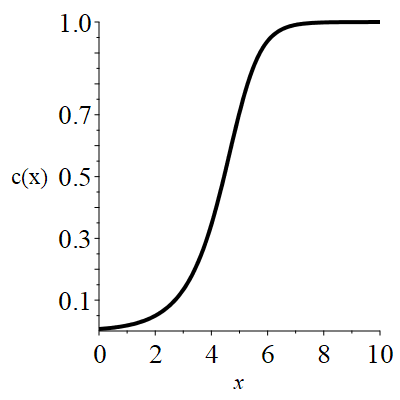}{a}
\includegraphics[width=4cm]{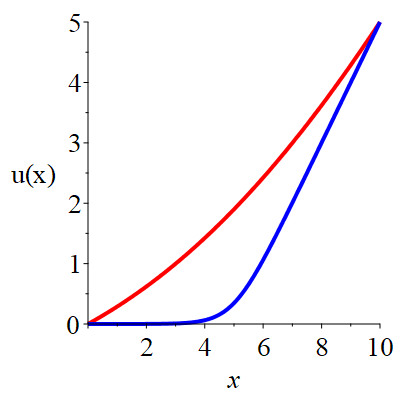}{b}
\caption{Fast path, $k=1$, and $x_0=5$.    a) The speed is a logistic function centered on $x_0$ and having the slope $\beta = 2$. b) The curves  u(x) when the slopes are $\beta=0.2$ (Red), $\beta=2$ (blue).} 
\label{FastPath}
\end{figure}

We consider another problem related to the Snell's law governing the refraction of light in its passage from one medium to another, provided that the observed refractive index of the medium is identified with the inverse of  speed $c^{-1}(x)$. Pierre de Fermat observed that Snell's law stems from the principle that light travels the path that takes the least amount of time. Specifically, the functional $F(x, u, u^{'})$ is the time taken to cover the path $u(x)$ for $x\in [a,b]$, i.e $F(x,u,u') = \frac{\sqrt{1+u'^2}}{c(x)}$ and $c(x)>0$. The  terms involved  in Eq.~\ref{odeFinal} are $\frac{ \partial F(x, u, u^{'})}{\partial u}  =0$, $\frac{ \partial F(x, u, u^{'})}{\partial u^{'}} = \frac{u'}{c^2(x) F(x, u, u^{'})}$,   $ \frac{d F(x, u, u^{'})}{d x}  = \frac{ u' u''}{c^2(x) F(x, u, u^{'})} -\frac{c'(x)}{c(x)} F(x, u, u^{'})$, and $ \frac{ d }{d x}  \frac{ \partial F(x, u, u^{'})}{\partial u^{'}} =\frac{-u'(x)c'(x)}{c^3(x) F(x, u, u^{'})}+\frac{u''(x)}{c^4(x)F^3(x, u, u^{'})}$.
Eq.~\ref{odeFinal} is rewritten as:
\begin{equation}
F(x, u, u^{'})( \frac{-u'(x)c'(x)}{c^3(x) F(x, u, u^{'})}+\frac{u''(x)}{c^4(x)F^3(x, u, u^{'})})-(\alpha-1) \frac{u'}{c^2(x) F(x, u, u^{'})} (\frac{ u' u''}{c^2(x) F(x, u, u^{'})} -\frac{c'(x)}{c(x)} F(x, u, u^{'}))=0
\label{odeFinala22}
\end{equation}
Eq.~\ref{odeFinala22} can be rewritten as $\frac{ \partial F(x, u, u^{'})}{\partial u^{'}}= k$, where $k$  is a constant when $\alpha=1$ because $u(x)$ is an ignorable coordinate. It follows that  the derivative  of the stationary curve  is $u'(x)= \mp kc(x)/ \sqrt{1- k^2 c^2(x)}$. For example, consider $c(x)=x$ and $k \neq 0$, then $u(x)=\mp \sqrt{1-k^2 x^2}/k+p$, which is the circle $(u(x)-p)^2 +  x^2 = k^{-2}$. This is well-known solution in the calculus of variations. When $\alpha=0$, Eq.~\ref{odeFinala00} is rewritten as $u''(x)(1- u'^2(x))=0$. The case $u'^2(x) =  1$ is a particular case of the $u''(x)=0$, the stationary curve $u(x)$ solution to the last ODE is independent of the speed;  $u(x) = s x +d$. The second variation in Eq.\ref{TestFinal0}, $\frac{\partial^2 F}{\partial u^2}=\frac{\partial^2 F}{\partial u \partial u'}=\frac{\partial F}{\partial u}=0$,  is   $(1-s^2) \int_a^b h'^2(x) dx/\delta(s^2+1)^2$. The second variation is positive when $s^2 < 1$  because $h(x)$ cannot be a constant, and therefore the stationary curve is a minimum in this case.   In other words, only straight lines with slopes in $]-1,1[$ are minimums. 
When $\alpha=2$, we can rewrite the ODE in Eq.~\ref{odeFinala22} as $\frac{d}{dx} F(x, u, u^{'}) \frac{ \partial F(x, u, u^{'})}{\partial u'}=0$; that is $F(x, u, u^{'}) \frac{ \partial F(x, u, u^{'})}{\partial u'}=k$ and therefore $u'(x)=kc^2(x)$. Let us set $c(x)=x$, then $u(x)=k x^3/3+p$. It is a minimum because the second variation  test in Eq.~\ref{TestFinalN0} leads to   $ 2 \int_{a}^b \frac{h'^2(x)}{x^2}/\delta dx > 0$.  If we consider that there are two mediums, having different refractive indexes whose variation  is smooth along the x-axis. The speed function $c(x)=\sqrt{e^{\beta (x-x_0)}/(1+e^{\beta (x-x_0)})}$, where  $\beta$ and $x_0$ are parameters describes two mediums.  By considering this speed function,  the stationary curve   $u(x)=k ln(1+e^{\beta (x-x_0)})/\beta+p$ is a minimum because the second variation test in Eq.~\ref{TestFinalN0} is rewritten as $2 \int_{a}^b h'^2(x) \sqrt{ \frac{1+e^{\beta (x-x_0)}}{e^{\beta (x-x_0)}}} dx>0$, where $w$ is a positive function of $a$, $b$, and $k^2$. This curve $u(x)$ is plotted in Fig.~\ref{FastPath}.

\section{Conclusion}
In this document, the calculus of variations is revisited from a new angle. Instead of the usual variational formulation, functional centrality based on H\"older's mean is used. For a given problem, the solution is a family of curves, where each one corresponds to a certain data selection criterion induced by the centrality measure used. Some of these solutions cannot be obtained with the usual variational formulation.
To illustrate, the fastest and shortest paths are studied.

\bibliographystyle{elsarticle-num}


\end{document}